\documentclass[11pt]{elsarticle}
\usepackage{amsmath}
\usepackage{amssymb}
\usepackage{graphicx}
\usepackage{amscd}
\usepackage{amsfonts}
\usepackage{bm}                                                
\usepackage{float}
\usepackage{latexsym}
\usepackage{lscape}
\usepackage{mathrsfs}
\usepackage{verbatim}
\usepackage{color}

\newcommand{\om}{[0,L]}
\newcommand{\ds}{\displaystyle}
\newtheorem{theorem}{Theorem}[section]

\newtheorem{lemma}{Lemma}[section]

\newcommand{\xx}{{x}}

\newcommand{\p}{\partial}

\date{\today}

\begin{document}

\begin{frontmatter}

\title{Uniqueness of determining the variable fractional order in variable-order time-fractional diffusion equations }

\author[1]{Xiangcheng Zheng}
\ead{xz3@math.sc.edu}
\author[2]{Jin Cheng}
\ead{jcheng@fudan.edu.cn}
\author[1]{Hong Wang}
\ead{hwang@math.sc.edu}

\address[1]{Department of Mathematics, University of South Carolina, Columbia, South Carolina 29208, USA}

\address[2]{School of Mathematical Sciences, Fudan University, 200433 Shanghai, China}

\begin{abstract}
We study an initial-boundary value problem of variable-order time-fractional diffusion equations in one space dimension. Based on the wellposedness of the proposed model and the smoothing properties of its solutions, which are shown to be determined by the behavior of the variable order at the initial time, a uniqueness result for an important inverse problem of determination of
the variable order in the time-fractional derivative contained in the proposed model from observations of its solutions is obtained. 
\end{abstract}


\end{frontmatter}


\section{Introduction}

Integer-order diffusion equations were derived under the assumptions that the underlying independent and identically distributed particle movements have (i) a mean free path and (ii) a mean waiting time \cite{MetKla00}, which hold for the diffusive transport of solutes in homogeneous aquifers when the solute plumes were observed to decay exponentially \cite{Bear61,Bear72}. However, field tests showed that the diffusive transport of solutes in heterogeneous aquifers often exhibit highly skewed power-law decaying tails \cite{BenSchMee,MeeSik,MetKla00}. Traditional practice to address the impact of the heterogeneity of the media is to tweak the variable parameters that multiply the pre-set integer-order diffusion equations to fit the training data, which tends to recover a rapidly varying, scale-dependent diffusivity and may overfit the training data but yield less accurate prediction on testing data \cite{SchBenMee03}. The reason is that the solutions of integer-order diffusion equations are Gaussian type and so catch exactly the exponential decaying behavior of the solute transport in homogeneous media but not the highly skewed power-law decaying behavior of the solute transport in heterogeneous media. 

Fractional diffusion equations (FDEs) were derived so that their solutions can have highly skewed power-law decaying tails and so catch the same behavior of the solute transport in heterogeneous media, which explains why FDEs model anomalously diffusive transport in heterogeneous aquifers more accurately than integer-order diffusion equations do \cite{MeeSik,MetKla00}. However, FDEs introduce new mathematical issues that are not common in the context of integer-order diffusion equations. Stynes et al. \cite{StyOriGra} proved that the first-order time derivatives of the solutions to the time-fractional diffusion equations (tFDEs) of order $0 < \alpha < 1$ have a singularity of order $O(t^{\alpha-1})$ near time $t=0$, which makes the error estimates in the literature that were proved under full regularity assumptions of the true solutions inappropriate. 

Nevertheless, the singularity of the solutions to the tFDEs at $t=0$ does not seem to be physically relevant to the subdiffusive transport they model. The current consensus why this phenomenon occurs lies in the incompatibility between the nonlocality of the power law decaying tails and the locality of the classical initial condition at the time $t=0$. However, the key issue of developing a physically relevant tFDE model that can correct the nonphysical behavior of the existing tFDEs remains unresolved. 

It was speculated in \cite{JMAA} that to eliminate their nonphysical singularity at the time $t=0$, as $t \rightarrow 0^+$ the power law decaying tail of the solutions to tFDEs should switch smoothly to an exponentially decaying tail to account for the impact of locality of the initial condition at the time $t=0$. That is, a physically relevant tFDE model should bear a variable-order near the time $t=0$, since the power of the power law decaying tail is determined by the order of the tFDEs. Moreover, variable-order tFDEs themselves occur in many applications \cite{LiWan,LorHar,SunChaZha,SunCheChe,ZenZhaKar,ZhuLiu}, as the order of tFDEs is closely related to the fractal dimension of the porous media determined by the Hurst index \cite{MeeSik} that changes as the geometrical structure or property of the media changes, e.g., in such applications as shale gas production due to the hydraulic fracturing and shape memory polymer due to the change of its miscrostructure \cite{LiWan,SunChaZha,SunCheChe}. 

Although variable-order tFDEs have appeared in increasingly more applications, their rigorous mathematical analysis is meager. In \cite{UmaSte}, a piecewise-constant order fractional ordinary differential equation was solved analytically on each time interval of a constant order via the Mittag-Leffler expression of the solution to the constant-order fractional differential equations \cite{Die}, with the solutions on the previous pieces and the solution value at the left end of the current piece as the source term and the initial data, respectively. Kian et al.  \cite{KiaSocYam} studied the wellposedness of a variable-order tFDE, in which the variable order is a function of the space variable, so the Laplace transform technique, which is a widely used analysis technique in the study of tFDEs, can be fully utilized to derive the desired results. In the variable-order tFDEs that occur in applications the variable order tends to be time dependent, hence the previous analytical techniques do not apply in the current context. New techniques have to be developed to analyze the wellposedness of these problems, and, probably more importantly, the regularity of the solutions to these problems.

In \cite{JMAA} we proved the wellposedness of the classical initial-boundary value problems of variable-order linear tFDEs in multiple space dimensions. We further proved that the regularity of their solutions depends on the behavior of the variable order (and its derivatives) at time $t=0$, in addition to the usual smoothness assumptions. More precisely, their solutions have the full regularity as their integer-order analogues do if the variable order has an integer limit at $t=0$ or have certain singularity at $t=0$ as the constant-order tFDEs do if the variable order has a non-integer value at time $t=0$. 

In most applications, the parameters in the governing diffusion equations, such as the diffusivity coefficients, source term, the boundary and initial data, are not given a priori. Rather, they have to be inferred from the measurements as an inverse problem by solving the forward problems repeatedly \cite{Kir}. In recent years the inverse problems of determining the parameters in (constant-order) tFDEs, in particular the order(s) of tFDEs that was not encountered in the inverse problems of integer-order diffusion equations, have attracted increasingly more research actitivities \cite{Jin12,KiaOksSoc,LiLiuYam19,LiuYam,XuCheYam,ZheWei}. In \cite{Che09} the uniqueness of the inverse problem of simultaneously determining the fractional order and the diffusivity coefficients of the one-dimensional tFDEs with the homogeneous right-hand side and variable diffusivity coefficients with the 
the observations at the left end point of the spatial interval, was proved. The inverse problem of determining the fractional order and a time varying kernel in tFDEs was studied in \cite{Jan}. A simultaneous inversion of the space-dependent diffusivity coefficient and the fractional order in one-dimensional tFDEs based on observations from the end points of the spatial interval was developed in \cite{LiZhaJia}. A numerical inversion of the fractional orders of the multi-term tFDEs in multiple space dimensions was studied in \cite{SunLiJia}.

To our best knowledge, up to now there is no mathematically proved result on the determination of the variable fractional order in variable-order tFDEs. The only available result is a numerical study of simultaneously determining the fractional order and the diffusion coefficient in a variable-order tFDE in which the variable order depends on both space and time variables \cite{WanWanLi}. In this paper we prove the uniqueness of determining the variable fractional order of the initial-boundary value problems of variable-order tFDEs in one space dimension with some available observed values of the unknown solutions inside the spatial domain. 

The rest of the paper is organized as follows: In Section 2 we discuss the modeling issues and present a physically relevant initial-boundary value problem of variable-order tFDEs. In Section 3 we address the wellposedness of the proposed model and the smoothing properties of its solutions, based on which we prove a uniqueness result in an inverse problem of determination of variable orders in the proposed model from some available observed values of the unknown solutions inside the spatial domain in Section 4.

\section{Modeling issues by tFDEs}

In this section we address the modeling issues of tFDEs. We begin with the widely used conventional tFDE model of order $0 < \alpha < 1$ \cite{MeeSik,MetKla00}
\begin{equation}\label{tFDE1}
{}_0^C D_t^{\alpha}u -  K\,u_{xx} = 0,
\end{equation}
where the Caputo fractional differential operator ${}_0^CD_t^{\alpha}$ is defined by \cite{SamKil}
\begin{equation}\label{tFDE1:e1}\begin{array}{rl}
\ds {}_0I_t^{\alpha}g(t) & \ds := \frac{1}{\Gamma(\alpha)} \int_0^t \frac{g(s)}{(t-s)^{1-\alpha}}ds, \\[0.125in]
\ds {}_0^C D_t^\alpha g(t) & \ds := {}_0I_t^{1-\alpha}g'(t) = \frac{1}{\Gamma(1-\alpha)} \int_0^t \frac{g'(s)}{(t-s)^{\alpha}}ds.
\end{array}\end{equation}
We note that the tFDE \eqref{tFDE1} was derived via a stochastic framework of the continuous time random walk formulation as the number of particle jumps tends to infinity (while the mean jump length shrinks) \cite{MeeSik,MetKla00}. In other words, the tFDE \eqref{tFDE1}, as the diffusion limit of the continuous time random walk in the phase space, holds for relatively large time $t > 0$, rather than all the way up to the time $t=0$ as often assumed in the literature. This gap explains partially why the solutions to the conventional initial-boundary value problem of the tFDE \eqref{tFDE1} exhibit nonphysical singularity near the initial time $t = 0$.

A modified tFDE model was derived in \cite{SchBenMee03} to model the anomalously diffusive transport of solute in heterogeneous porous media
\begin{equation}\label{tFDE2}\begin{array}{c}
u_t + k\;{}_0^C D_t^{\alpha} u -  K\, u _{xx} = 0,~~(x,t) \in [0,L]\times(0,T]; \\[0.05in]
\ds u(x,0)=u_0(x),~x\in [0,L],\quad u(0,t)=u(L,t) = 0,~t \in [0,T]. 
\end{array}\end{equation}
Here $u_t$ is the first-order derivatives of $u$ and $K>0$ is the diffusivity coefficient. 

With the prescribed initial distribution $u_0(\xx)$ of the solute that is often present in the water phase initially, the governing tFDE \eqref{tFDE2} models the anomalously diffusive transport of the solute in the heterogeneous porous media. During the transport process, a large amount of solute may get absorbed to the solid matrix due to the impact of adsorption that deviates from the transport of the solute in the bulk fluid phase. The absorbed solute may get slowly released later to the bulk fluid phase that leads to a subdiffusive transport process, which makes the remediation effort of the contaminant solute in groundwater much less effective than it is observed in laboratory experiments. In fact, some studies suggested that the remediation of contaminated aquifers may take decades or centuries longer than integer-order diffusione equation models had predicted  \cite{MeeSik}. 

In the traditional integer-order diffusion equation models, the amount of the adsorbed solute was expressed as a (linear or nonlinear) function of the amount of the solute in the bulk phase to account for the impact of adsorption and desorption. This would yield a retardation coefficient in front of the first-order time derivative term in the integer-order diffusion equation \cite{Bear72}. As it relates only the amount of adsorbed solute in the solid matrix to the amount of solute in the bulk phase at the current time instant, the conventional integer-order diffusion equation model does not account for the amount accumulated to the solid matrix. In the tFDE model \eqref{tFDE2}, the $u_t$ term models the Fickian diffusive transport of the solute in the bulk fluid phase, which consists of $1/(1+k)$ portion of the total solute mass, and the $k\;{}_0^C D_t^{\alpha} u$ term models the subdiffusive transport of the solute absorbed to the solid matrix, which is $k/(1+k)$ portion of the total solute mass. Note that the governing tFDE \eqref{tFDE2} holds on the entire time interval including the initial time $t=0$.

A possible remedy to eliminate the nonphysical singularity of the solutions to the initial-boundary value problem of the tFDE \eqref{tFDE1} (and \eqref{tFDE2}) as proved in \cite{StyOriGra} is to vary their power law decaying tails smoothly to exponentially decaying tails as $t \rightarrow 0^+$ to account for the impact of locality of the initial condition at the time $t=0$. This leads to variable-order tFDEs. 

In shale gas production, shale formation often has insufficient permeability due to the existence of fine-scale pores that results in a large amount of shale gas molecules absorbed to the surface of the pore-throat structures and significant decrease of flow rate to the wellbore. Hydraulic fracturing is often used to increase the pore sizes and so the kermeability of the shale formation and to increase shale gas production, which changes the fractal dimension of the shale formation and so the order of the tFDEs \cite{LeeBoc,MeeSik}. This again leads to variable-order tFDEs.

Motivated by these observations, in this paper we consider the initial-boundary value problem of the following variable-order linear tFDE 
\begin{equation}\label{Model}\begin{array}{c}
u_t + k(t)\; {}_0^C D_t^{\alpha(t)} u -K \,u_{xx}= 0,~~(x,t) \in [0,L]\times(0,T]; \\[0.1in]
u(x,0)=u_0(x),~x\in [0,L], \quad u(0,t)=u(L,t) = 0,~t \in [0,T]. 
\end{array}\end{equation}
Here ${}_0I_t^{\alpha(t)}$ and ${}_0^C D_t^{\alpha(t)}$ denote the variable-order fractional integral and Caputo fractional differential operators, respectively \cite{LorHar,SunChaZha,SunCheChe}
\begin{equation*}\begin{array}{rl}
\ds {}_0 I_t^{\alpha(t)}g(t) &\ds := \frac{1}{\Gamma(\alpha(t))}\int_0^t\frac{g(s)}{(t-s)^{1-\alpha(t)}}ds,\\[0.125in]
\ds {}_0^C D_t^{\alpha(t)}g(t) & \ds :={}_0 I_t^{1-\alpha(t)}g'(t) = \frac{1}{\Gamma(1-\alpha(t))}\int_0^t\frac{g'(s)}{(t-s)^{\alpha(t)}}ds,
\end{array}\end{equation*}
which are a variable-order analogue of the constant-order fractional integral and differential operators defined in \eqref{tFDE1:e1}.

\section{Wellposedness and smoothing properties of the variable-order linear tFDE \eqref{Model}}

Let $m \in \mathbb{N}$, the set of nonnegative integers and let ${\cal I} \subset \mathbb{R}$ be a bounded (open or closed or half open and half closed) interval. Let $C^m(\cal I)$ be the spaces of continuous functions with continuous derivatives up to order $m$ on $\cal I$ and $C(\mathcal I) := C^0(\mathcal I)$, equipped with the standard norms. Let $H^{m}(0,L)$ be the Sobolev spaces of Lebesgue square integrable functions with their weak derivatives of order $m$ being in $L^2(0,L)$. Let $H^m_0(0,L)$ be the completion of $C^\infty_0(0,L)$, the space of infinitely differentiable functions with compact support in $\om$, in $H^m(0,L)$. For non-integer $s\geq 0$, the fractional Sobolev spaces $H^s(0,L)$ are defined by interpolation. All the spaces are equipped with the standard norms \cite{AdaFou}. 

For a Banach space ${\cal X}$ the norm $\| \cdot \|_{\cal X}$, let $C^m(\mathcal I;\mathcal X)$ be the Banach spaces equipped with the norm unctions on the interval $\cal I$ \cite{AdaFou,Eva} 
$$\begin{array}{c}
\ds C^m(\mathcal I;\mathcal X) := \bigg \{w(x,t) : \Bigl \| \frac{\p^l w}{\p t^l} \Bigr \|_{\mathcal X} \in C^l({\cal I}), \quad l = 0,1,\ldots,m \bigg\},\\[0.15in]
\ds \| w \|_{C^m({\mathcal I};\mathcal X)} := \max_{0 \le l \le m} \sup_{t \in {\mathcal I}} ~ \Bigl \| \frac{\p^l w}{\p t^l} \Bigr \|_{\mathcal X}. 
\end{array}$$
To better characterize the temporal singularity of the solution at the initial time, we define the weighted Banach spaces involving time $C^m_\gamma((0,T];\cal X)$ with $m\geq 2$, $0\leq \mu <1$ modified from those in \cite{Mcl}
$$\begin{array}{c}
\ds C^m_\mu((0,T];\mathcal X):= \big\{ w \in C^1([0,T];\mathcal X) : \| w \|_{C_\mu^m((0,T];\mathcal X)}<\infty \big\},\\[0.075in]
\ds \| w \|_{C_\mu^m((0,T];\mathcal X)} := \| w \|_{C^{1}([0,T];\mathcal X)} + \sum_{l=2}^{m}  \sup_{t \in (0,T]} t^{l-1-\mu}\Bigl \| \frac{\p^l w}{\p t^l} \Bigr \|_{\mathcal X}.
\end{array}$$

Let $\{\lambda_i,\phi_i \}_{i=1}^\infty$ be the eigenvalues and eigenfunctions of the Sturm-Liouville problem 
\begin{equation}\label{SturmLiouville}
-  K\, D_{xx}\phi_{i}(\xx) = \lambda_i \phi_i(\xx), ~\xx \in [0,L]; \quad \phi_i(0)=\phi_i(L) = 0. 
\end{equation}
It is known that \cite{Eva}
\begin{equation}
\ds \lambda_i=\frac{K\,i^2\pi^2}{L^2},~~\phi_i(x)=\sin\big(\frac{i\pi x}{L}\big),~~i\in \mathbb N^+.
\end{equation}
 Furthermore, for any $\gamma \ge 0$ the Sobolev space defined by 
$$\begin{array}{c}
\ds \hspace{-0.1in} \check{H}^{\gamma}(0,L) := \Big \{v \in L^2(0,L): | v |_{\check{H}^\gamma}^2 := \big ((-D_{xx})^\gamma v,v\big) = \sum_{i=1}^{\infty} \lambda_i^{\gamma} (v,\phi_i)^2 < \infty \Big \}
\end{array}$$ 
is a subspace of the fractional Sobolev space $H^\gamma(0,L)$ characterized by \cite{AdaFou,SakYam,Tho}
$$
\check{H}^\gamma(0,L) = \big \{v \in H^\gamma(0,L): ~(-D_{xx})^l v(0)= (-D_{xx})^l v(L)= 0,~~ l < \gamma/2 ,~~l \in \mathbb N \big \}
$$
and the seminorms $| v |_{\check{H}^\gamma}$ and $| v |_{{H}^\gamma}$ are equivalent in $\check{H}^\gamma$.

We make the following assumptions on the data:
\paragraph{\textbf{Assumption A}}
Suppose that $\alpha, k \in C[0,T]$ and
\begin{equation}\label{alpha}
0\leq \alpha(t)\leq \alpha_*<1,~t\in [0,T],~~\lim_{t \rightarrow 0^+} \big ( \alpha(t) - \alpha(0) \big) \ln t ~\mbox{exists}.
\end{equation}

Then the following theorems hold \cite{JMAA}.

\begin{theorem}\label{thm:Wellpose}
Suppose that Assumption A holds, and $u_0 \in \check{H}^{\gamma + 2}$ for some $1/2 < \gamma \in \mathbb{R}^+$. Then problem (\ref{Model}) has a unique solution $u \in C^1\big([0,T];H^\gamma \big)$ and the following stability estimates hold \begin{equation*}
\begin{array}{c}
\|u\|_{C([0,T];H^\gamma(0,L))} \leq Q\|u_0\|_{\check H^\gamma(0,L)}, \quad
\|u\|_{C^1([0,T];H^\gamma(0,L))}  \le Q \|u_0\|_{\check H^{2+\gamma}(0,L)}.
\end{array}
\end{equation*}
\end{theorem}

\begin{theorem}\label{thmreg1}
For $2 \le n \in \mathbb{N}$, suppose that $u_0\in \check{H}^{\gamma + 2n}$, $\alpha , k \in C^{n-1}[0,T]$, and (\ref{alpha}) holds. 

If $\alpha^{(l)}(0)=0$ for $0 \leq l \leq n-2$ and $\lim_{t \rightarrow 0}\alpha^{(n-2)}(t)\ln t$ exists, then $u \in C^n([0,T];\check H^\gamma(0,L))$ 
such that 
$$\|u\|_{C^n([0,T];\check H^\gamma(0,L)} \le Q \|u_0\|_{\check H^{\gamma+2n}(0,L)}.$$ 
 
If $\alpha(0) > 0$, then $u\in C^n((0,T];\check H^\gamma(0,L)) \cap C^{n}_{1-\alpha(0)}((0,T];\check H^\gamma(0,L))$ and
$$\|u\|_{C^{n}_{1-\alpha(0)}((0,T];\check H^\gamma(0,L))} \le Q \|u_0\|_{\check H^{\gamma+2n}(0,L)}.$$ 
\end{theorem}

\section{The inverse problem of determining the variable order in variable-order tFDEs}

In this section we prove the main result of this paper, the uniqueness of the inverse problem of determining the variable order in the variable-order tFDE model in (\ref{Model}) based on some observations of the solution $u(x,t)$ to the initial-boundary value problem \eqref{Model}.
We first prove a lemma for future use.
\begin{lemma}\label{Lem}
Assume that 
\begin{equation}\label{observe}
\ds \sum_{i=1}^\infty g_i(t)\phi_i(x)=0,~~(x,t)\in (a,b)\times [0,T]
\end{equation}
 for some $(a,b)\subset [0,L]$. Then $g_i(t)\equiv 0$ on $t\in [0,T]$ for $i\in \mathbb N^+$.
\end{lemma}
\textbf{Proof.}
 Let $\bm G:=(g_i(t)\phi_i(x))_{i=1}^\infty$. Then the relation (\ref{observe}) can be rewritten in the vector form by $\bm \Lambda_0^\top \bm G=0$ where $\bm \Lambda_0=(\cdots,1,\cdots)$. We then apply $K\,D_{xx}$ on both sides of (\ref{observe}) and use (\ref{SturmLiouville}) to obtain
 \begin{equation}\label{o1}
\ds \sum_{i=1}^\infty g_i(t)\lambda_i\phi_i(x)=0,~~(x,t)\in (a,b)\times [0,T],
\end{equation}
which can be rewritten as $\bm \Lambda_1^\top \bm G=0$ with $\bm \Lambda_1=(\lambda_i)_{i=1}^\infty$. After repeating this procedure for $n$ times we have $\bm \Lambda_n^\top \bm G=0$ with $\bm \Lambda_n=(\lambda_i^n)_{i=1}^\infty$. We collect the relations $\bm \Lambda_n^\top \bm G=0$ for $n\in \mathbb N$ to get an infinite dimension linear system with a Vandermonde coefficient matrix 
\begin{equation}\label{o2}
\bm\Lambda\,\bm G=\bm 0,~~\bm\Lambda:=(\bm\Lambda_n^\top)_{n=0}^\infty.
\end{equation}

As $\lambda_i\neq \lambda_j$ for $i\neq j$, (\ref{o2}) yields $\bm G\equiv \bm 0$, i.e., $g_i(t)\phi_i(x)=0$ on $(x,t)\in (a,b)\times [0,T]$ for $i\in \mathbb N^+$. For any $i\in \mathbb N^+$, there exists an $x_i\in (a,b)$ such that $\phi_i(x_i)\neq 0$, which immediately leads to $g_i(t)\equiv 0$ on $t\in [0,T]$. So the proof is finished.
$\blacksquare$

Let the admissible set $\mathcal A$ be defined by
\begin{equation*}
\mathcal A:=\big\{\alpha(t):\alpha(t)~\mbox{is analytic on}~[0,T]~\mbox{and satisfies} ~\eqref{alpha} \big\}. 
\end{equation*}
Many types of functions can be taken to fulfill the constraints of the admissible set $\mathcal A$, e.g., the polynomials that satisfy \eqref{alpha}.

We prove the main result of this work in the following theorem.
\begin{theorem}\label{thm:inverse}
Suppose that Assumptions A holds, $k(0)\neq 0$, and $0 \not \equiv u_0 \in H^{\gamma+2}(0,L)$ for some $\gamma > 1/2$. Then the variable order $\alpha \in \mathcal A$ in the initial-boundary value problem \eqref{Model} is determined uniquely.  Namely, let $\hat \alpha \in \mathcal A$ and $\hat u(x,t)$ be the solution to the problem  
\begin{equation}\label{Model1}\begin{array}{c}
\hat u_t + k(t)\; {}_0^C D_t^{\hat \alpha(t)} \hat u - K\,\hat u_{xx}= 0,~~(x,t) \in [0,L]\times(0,T]; \\[0.1in]
\hat u(x,0)=u_0(x),~x\in [0,L]; \quad \hat u(0,t)=\hat u(L,t) = 0,~t\in [0,T]. 
\end{array}\end{equation}
If
\begin{equation}\label{Inv:e3}
u(x,t) = \hat u(x,t), \quad (x,t) \in (a,b)\times [0,T]
\end{equation}
for some $(a,b)\subset [0,L]$, then we have
\begin{equation*}
\alpha(t) = \hat \alpha(t), \quad t \in [0,T].
\end{equation*} 
\end{theorem}

\textbf{Proof.} We express $\hat u$ and the solution $u$ to problem \eqref{Model} in terms of $\{\phi_i\}_{i=1}^\infty$ as follows \cite{Luc,SakYam,StyOriGra} 
\begin{equation}\label{Inv:e5}\begin{array}{l}
\ds u(x,t)=\sum_{i=1}^\infty u_i(t)\phi_i(x), ~~ u_i(t) := \big (u(\cdot,t),\phi_i \big),~~t \in [0,T];\\[0.2in]
\ds \hat u(x,t)=\sum_{i=1}^\infty \hat u_i(t)\phi_i(x),~~\hat u_i(t) := \big (\hat u(\cdot,t),\phi_i \big),~~t \in [0,T].
\end{array}\end{equation}

We plug these expansions into (\ref{Model}) and use (\ref{SturmLiouville}) to obtain 
\begin{equation*}
\begin{array}{l}
\ds \sum_{i=1}^{\infty} u_i'(t) \phi_i(x) + k(t) \sum_{i=1}^{\infty} {}_0D_t^{\alpha(t)} u_i(t)\phi_i(x) = -\sum_{i=1}^{\infty}\lambda_i u_i(t) \phi_i(x), \\[0.175in]
\ds \quad\quad \quad \forall x \in \om, ~t \in (0,T].
\end{array}
\end{equation*}
Hence, $u$ is the solution to problem (\ref{Model}) if and only if $\{u_i\}_{i=1}^\infty$ satisfy
\begin{equation}\label{Inv:e7}\begin{array}{c}
\ds u_i'(t)+k(t) {}_0D_t^{\alpha(t)} u_i(t) = -\lambda_i u_i(t), ~~t \in (0,T],\\[0.1in]
\ds  u_i(0) = u_{0,i} := (u_0,\phi_i), ~~i =1,2,\cdots.
\end{array}\end{equation}
Similarly, $\hat u$ is the solution to problem \eqref{Model1} if and only if $\hat u_i(t)$ satisfy  
\begin{equation}\label{Inv:e8}\begin{array}{c}
\ds \hat u_i'(t)+k(t) {}_0D_t^{\hat\alpha(t)} \hat u_i(t) = -\lambda_i\hat u_i(t), ~~t \in (0,T],\\[0.1in]
\ds \hat u_i(0) = u_{0,i}, ~~i=1,2,\cdots.
\end{array}\end{equation}

By (\ref{Inv:e3}) and (\ref{Inv:e5}) we have
\begin{equation*}
\begin{array}{c}
\ds 0=u(x,t)-\hat u(x,t)=\sum_{i=1}^\infty (u_{i}(t)-\hat u_{i}(t))\phi_{i}(x),  \quad  (x,t) \in (a,b)\times [0,T],
\end{array}
\end{equation*}
which, by Lemma \ref{Lem}, implies that $u_{i}(t)=\hat u_{i}(t)$ on $t\in [0,T]$ for any $i\in \mathbb N^+$.

As $u_0(x)\not\equiv 0$, there exists at least one $i_*\in \mathbb N^+$ such that $u_{0,i_*}\neq 0$.  Then we replace $\hat u_i(t)$ by $u_i(t)$ for $i=i_*$ in the first equation of (\ref{Inv:e8}) and then subtract it from the first equation in (\ref{Inv:e7}) to obtain
\begin{equation*}
k(t)\big({}^C_0D^{\alpha(t)}_t-{}_0^CD^{\hat\alpha(t)}_t\big)u_{i_*}(t)=0,~~t\in (0,T],~~u_{i_*}(0)=u_{0,i_*}.
\end{equation*}
By the assumptions $k(0)\neq 0$ and $k(t)\in C[0,T]$, there exists an $0<\varepsilon_1<1$ such that 
\begin{equation}\label{thminv:e3}
\big({}_0^CD^{\alpha(t)}_t-{}_0^CD^{\hat\alpha(t)}_t\big)u_{i_*}(t)=0,~~t\in (0,\varepsilon_1],~~u_{i_*}(0)=u_{0,i_*}.
\end{equation}
For any fixed $t\in (0,\varepsilon_1]$, we consider the variable-order fractional derivative ${}_0^CD_t^{\alpha(t)}u_{i_*}(t)$ as a function of $\alpha(t)$ and then apply the following relation
$$\begin{array}{l}
\ds \frac{d}{d\alpha(t)}\big({}_0^CD_t^{\alpha(t)}u_{i_*}(t)\big)=\frac{d}{d\alpha(t)}\Big(\frac{1}{\Gamma(1-\alpha(t))}\int_0^t\frac{u_{i_*}'(s)ds}{(t-s)^{\alpha(t)}}\Big)\\[0.15in]
\qquad\ds=\int_0^t\Big(\frac{\Gamma'(1-\alpha(t))}{\Gamma^2(1-\alpha(t))}\frac{1}{(t-s)^{\alpha(t)}}-\frac{1}{\Gamma(1-\alpha(t))}\frac{\ln(t-s)}{(t-s)^{\alpha(t)}}\Big)u'_{i_*}(s)ds\\[0.15in]
\ds\qquad=\frac{1}{\Gamma(1-\alpha(t))}\int_0^t\Big(\frac{\Gamma'(1-\alpha(t))}{\Gamma(1-\alpha(t))}-\ln(t-s)\Big)\frac{u'_{i_*}(s)}{(t-s)^{\alpha(t)}}ds
\end{array}$$
and the mid-point formula on (\ref{thminv:e3}) to obtain
\begin{equation}\label{thminv:e4}
\begin{array}{l}
\ds  \Big[\frac{1}{\Gamma(1-\bar\alpha(t))}\int_0^t\Big(\frac{\Gamma'(1-\bar\alpha(t))}{\Gamma(1-\bar\alpha(t))}-\ln(t-s)\Big)\\[0.1in]
\ds\hspace{1.7in}\times\frac{u'_{i_*}(s)}{(t-s)^{\bar\alpha(t)}}ds\Big]\,(\alpha(t)-\hat\alpha(t))=0, 
\end{array}
\end{equation}
on $t\in (0,\varepsilon_1]$ where $\bar\alpha(t)$ lies in between $\alpha(t)$ and $\hat\alpha(t)$ for any $t\in (0,\varepsilon_1]$. 

As $u_{0,i_*}\neq 0$, we assume $u_{0,i_*}>0$ without loss of generality. By Theorem \ref{thm:Wellpose}, $u(x,t)$ and $u_t(x,t)$ are continuous in time, which implies that $u_{i_*}(t)=(u(\cdot,t),\phi_{i_*})$ and $u_{i_*}'(t)=(u_t(\cdot,t),\phi_{i_*})$ are also continuous in time on $t\in [0,T]$. Then there exists a positive $\varepsilon_2\leq \varepsilon_1$ such that $u_{i_*}(t)\geq \sigma>0$ on $[0,\varepsilon_2]$ for some $\sigma>0$. Furthermore, for the first equation of (\ref{Inv:e8}) with $i=i_*$, as $t$ tends to $0$, the second term on its left-hand side will tends to $0$ by the continuity of $u'_{i_*}(t)$ on $t\in [0,T]$ and hence the integrability of the kernel of integral. Then, based on this equation with $i=i_*$, there exists a positive $\varepsilon_3\leq \varepsilon_2$ such that 
\begin{equation}\label{thminv:e5}u'_{i_*}(t)\leq -\lambda_{i_*}\sigma/2<0,~~t\in(0,\varepsilon_3].\end{equation}

As $\alpha(t)$ and $\hat\alpha(t)$ are bounded away from $1$, $\bar\alpha(t)$ also has a positive upper bound less than $1$, which leads to the following estimate
\begin{equation*}
\ds \Big|\frac{\Gamma'(1-\bar\alpha(t))}{\Gamma(1-\bar\alpha(t))}\Big|\leq Q_0,~~t\in (0,\varepsilon_1].
\end{equation*} 
For this constant $Q_0$, there exists a positive $\varepsilon_4\leq \varepsilon_3$ such that $\ln t<-Q_0$ on $t\in (0,\varepsilon_4]$. which implies that 
$$\frac{\Gamma'(1-\bar\alpha(t))}{\Gamma(1-\alpha(t))}-\ln(t-s)>0,~~0<s<t,~~t\in (0,\varepsilon_4].$$
We incorporate this with (\ref{thminv:e5}) to conclude that the integral on the left-hand side of (\ref{thminv:e4}) is negative on $t\in (0,\varepsilon_4]$. Hence the only possible way that the equation (\ref{thminv:e4}) holds true is $\alpha(t)-\hat\alpha(t)=0$ on $t\in (0,\varepsilon_4]$. As $\alpha(t),\hat\alpha(t)\in \mathcal A$, the set of analytic functions satisfying (\ref{alpha}), we obtain $\alpha(t)=\hat\alpha(t)$ on $t\in [0,T]$, which finishes the proof.
$\blacksquare$

 \section*{Acknowledgements}
 This work was funded by the OSD/ARO MURI Grant W911NF-15-1-0562 and by the National Science Foundation under Grant DMS-1620194.

\end{document}